\documentclass{article}

\usepackage[utf8]{inputenc}
\usepackage[fleqn]{amsmath}
\usepackage{amsfonts, amssymb}

\usepackage[ntheorem]{empheq}
\usepackage[hyperref]{ntheorem}

\theoremstyle{hyphen}
\newtheorem{theorem}{Theorem}
\newtheorem{lemma}{Lemma}

\newtheorem{corollary}{Corollary}

\usepackage{extarrows} 
\usepackage{units} 
\usepackage{bm} 
\usepackage{bbold} 
\usepackage{tensor} 
\usepackage{soul} 
\usepackage{cancel} 
\usepackage[retainorgcmds]{IEEEtrantools} 
\IEEEeqnarraydefcolsep{0}{\leftmargini} 
\usepackage{graphicx}
\usepackage{fixltx2e} 

\usepackage{array} 
\newcolumntype{C}[1]{>{\centering\let\newline\\\arraybackslash\hspace{0pt}}m{#1}}

\usepackage[pdftex,colorlinks,bookmarks=false,pageanchor=false]{hyperref}
 \hypersetup{
     colorlinks,
     linkcolor={red!50!blue},
     citecolor={blue!50!green},
     urlcolor={blue!80!black}
 }
\usepackage[numbers,sort&compress]{natbib} 

\usepackage{color}
\usepackage{realboxes}
\usepackage{framed}
\usepackage{array} 
\usepackage{mathtools} 
\usepackage{mhchem} 

\usepackage{enumerate} 
\usepackage{enumitem}

\usepackage{bookmark} 
\usepackage{float} 
\usepackage{rotating}

\usepackage{ytableau} 
\usepackage{youngtab}
\usepackage{adjustbox} 

\ytableausetup{mathmode, centertableaux}
\usepackage{xargs}
\usepackage{letltxmacro} 
\usepackage{xcolor}
\usepackage{tikz}
\usepackage{tkz-graph} 
\usepackage{tkz-euclide}
\usetikzlibrary{matrix,fit,arrows,decorations}
\usepackage{pbox}

\newcommand{\qed}{\hfill\tikz{\draw[draw=black,line width=0.6pt] (0,0) rectangle (2.8mm,2.8mm);}\bigskip}

\usepackage[a4paper,left=2.0cm]{geometry}
\newcommand{\blockmatrix}[9]{
  \draw[draw=#4,fill=#5,every node/.style={inner sep=0,outer sep=0}] (0,0) rectangle( #1,#2);
  \ifthenelse{\equal{#6}{true}}
  {
    \draw[draw=#7,fill=#8] (0,#2) -- (#9,#2) -- ( #1,#9) -- ( #1,0) -- ( #1 - #9,0) -- (0,#2 -#9) -- cycle;
  }
  {}
  \draw ( #1/2, #2/2) node[inner sep=0,outer sep=0]{ #3};
}

\newcommand{\mblockmatrix}[4][none]{
  \begin{tikzpicture} 
  \ifthenelse{\equal{#1}{none}}
  {
    \blockmatrix{#2}{#3}{#4}{none}{none}{false}{none}{none}{0.0}
  }
  {
    \definecolor{fillcolor}{rgb}{#1}
    \blockmatrix{#2}{#3}{#4}{none}{fillcolor}{false}{none}{none}{0.0}
  }
  \end{tikzpicture}
}

\newcommand{\fblockmatrix}[4][none]{%
  \begin{tikzpicture}[outer sep=0,inner sep=0] 
  \ifthenelse{\equal{#1}{none}}{\blockmatrix{#2}{#3}{#4}{black}{none}{false}{none}{none}{0.0}}{\definecolor{fillcolor}{rgb}{#1}\blockmatrix{#2}{#3}{#4}{black}{fillcolor}{false}{none}{none}{0.0}}\end{tikzpicture}
}

\newcommand{\dblockmatrix}[4][none]{
  \begin{tikzpicture} 
  \ifthenelse{\equal{#1}{none}}
  {\blockmatrix{#2}{#3}{#4}{black}{none}{true}{black}{none}{0.35cm}}{\definecolor{fillcolor}{rgb}{#1}\blockmatrix{#2}{#3}{#4}{black}{none}{true}{black}{fillcolor}{0.35cm}}\end{tikzpicture}
}

\newcommand{\diagonalblockmatrix}[5][none]{
  \begin{tikzpicture} 

  \ifthenelse{\equal{#1}{none}}
  {
    \blockmatrix{#2}{#3}{#4}{black}{none}{true}{black}{none}{#5}
  }
  {
    \definecolor{fillcolor}{rgb}{#1}
    \blockmatrix{#2}{#3}{#4}{black}{none}{true}{black}{fillcolor}{#5}
  }

  \end{tikzpicture}
}

\renewenvironment{abstract}{%
\centering\begin{minipage}{.95\textwidth}
\sffamily{\bf Abstract:}}
{
\end{minipage}\vskip 3em}
\makeatletter
\renewcommand\@maketitle{%
\hfill
\begin{minipage}{\textwidth}
\vskip 2em
\let\footnote\thanks 
{\LARGE \bf \@title \par }
\vskip 1.5em
{\large \@author \par}
\end{minipage}
\vskip 3em \par
}
\makeatother
\usepackage{authblk}

\allowdisplaybreaks[1]

\title{A simple counting argument of the irreducible representations
  of $\mathsf{SU}(N)$ on mixed product spaces}
\author{J. Alcock-Zeilinger and H. Weigert}
\date{}

\begin{document}
\maketitle


\newcommand{\ud}{\mathrm{d}}
\newcommand{\bra}[1]{\langle #1 \vert}
\newcommand{\ket}[1]{\vert #1 \rangle}
\newcommand{\ScalarProd}[2]{\left\langle #1 \middle| #2 \right\rangle}
\newcommand{\Tr}[1]{\text{tr}\left(#1\right)}

\newcommand{\FPic}[2][{}]{\hspace{-0.27mm}\pbox{\textwidth}{\includegraphics[#1]{{#2}}}\hspace{-0.27mm}}

\newcommand{\ybox}[1]{\begin{ytableau} #1 \end{ytableau}}

\newcommand{\SUN}{\mathsf{SU}(N)}
\newcommand{\SU}[1]{\mathsf{SU}(#1)}
\newcommand{\GLN}{\mathsf{GL}(N)}
\newcommand{\MixedPow}[2]{V^{\otimes
    #1}\otimes\left(V^*\right)^{\otimes #2}}
\newcommand{\Pow}[1]{V^{\otimes #1}}
\newcommand{\DPow}[1]{\left(V^*\right)^{\otimes #1}}
\newcommand{\Lin}[1]{\mathrm{Lin}\left( #1 \right)}
\newcommand{\API}[1]{\mathsf{API}\left( #1 \right)}
\newcommand{\InvAlg}[1]{A\left[ S_{#1} \right]}
\newcommand{\Rsim}{\stackrel{\mathcal{R}}{\sim}}

\newcommand{\Ampc}[1]{\cancel{#1}_c\left[\mathcal{R}\right]}
\newcommand{\Ampr}[1]{\cancel{#1}_r\left[\mathcal{C}\right]}

\def\smath#1{\text{\scalebox{.8}{$#1$}}}
\def\sfrac#1#2{\smath{\frac{#1}{#2}}}

\def\Smath#1{\text{\scalebox{0.9}{$#1$}}}
\def\Sfrac#1#2{\Smath{\frac{#1}{#2}}}

\let\oldytableau\ytableau
\let\endoldytableau\endytableau
\renewenvironment{ytableau}{\begin{adjustbox}{scale=.78}\oldytableau}{\endoldytableau\end{adjustbox}}

\makeatletter
\newcommand{\vast}{\bBigg@{3}}
\newcommand{\Vast}{\bBigg@{4}}
\makeatother

\begin{abstract}
  That the number of irreducible representations of the special
  unitary group $\SUN$ on $\Pow{k}$ (which is also the number of Young
  tableaux with $k$ boxes) is given by the number of involutions in
  $S_k$ is a well known result (see, e.g.,~\cite{Knuth:1998} and other
  standard textbooks). In this paper, we present an alternative proof
  for this fact using a basis of projection and transition
  operators~\cite{Alcock-Zeilinger:2016sxc,Alcock-Zeilinger:2016cva}
  of the algebra of invariants of $\SUN$ on $\Pow{k}$. This proof is
  shown to easily generalize to the irreducible representations of
  $\SUN$ on mixed product spaces $\MixedPow{m}{n}$, implying that the
  number of irreducible representations of $\SUN$ on a product space
  $\MixedPow{m}{n}$ remains unchanged if one exchanges factors $V$ for
  $V^*$ and vice versa, as long as the total number of factors remains
  unchanged, \emph{c.f.} Corollary~\ref{cor:Irreps-SUN-MixedSpace}.
\end{abstract}


\section{Introduction}\label{sec:Introduction}

The aim of this paper is to provide an alternative, compact proof of
Theorem~\ref{thm:Number-Involutions-Irreps}, which states that the
number of irreducible representations of the special unitary group
$\SUN$ on a product space $\Pow{k}$ (for some positive integer $k$)
is the same as the number of self-inverse permutations
(\emph{involutions}) in the symmetric group $S_k$.; our proof is
presented in section~\ref{sec:Invariant-Theory}. We then go beyond
Theorem~\ref{thm:Number-Involutions-Irreps} and give a counting
argument for the irreducible representations of $\SUN$ on mixed
product spaces $\MixedPow{m}{n}$, \emph{c.f.}
Corollary~\ref{cor:Irreps-SUN-MixedSpace} in
section~\ref{sec:Irreps-MixedSpace}.

Before presenting our proof of
Theorem~\ref{thm:Number-Involutions-Irreps}, we outline, without
detail, the steps involved in the standard proof of this theorem (see,
e.g.,~\cite{Knuth:1998} and other standard textbooks). This standard
proof critically involves Young tableaux~\cite{Young:1928}, which are
combinatorial objects that classify the irreducible representations of
$\SUN$ on $\Pow{k}$. A Young tableau of size $k$ is an arrangement of
$k$ boxes that are left justified and top justified, and each box is
filled with a unique integer in $\lbrace1,2,\ldots,k\rbrace$ such that
that the numbers increase from left to right and from top to
bottom~\cite{Fulton:1997,Sagan:2000}; for example
\begin{equation}
  \begin{ytableau}
    1 & 3 & 4 & 9 \\
    2 & 7 \\
    5 & 8 \\
    6
  \end{ytableau}
\end{equation}
is a Young tableau. We denote the set of all Young tableaux consisting
of $k$ boxes by $\mathcal{Y}_k$. As already stated, the Young tableaux
in $\mathcal{Y}_k$ classify the irreducible representations of $\SUN$
on
$\Pow{k}$~\cite{Young:1928,Fulton:1997,Fulton:2004,Cvitanovic:2008zz,Tung:1985na}.

\begin{figure}[t]
  \centering
  \includegraphics[scale=0.85]{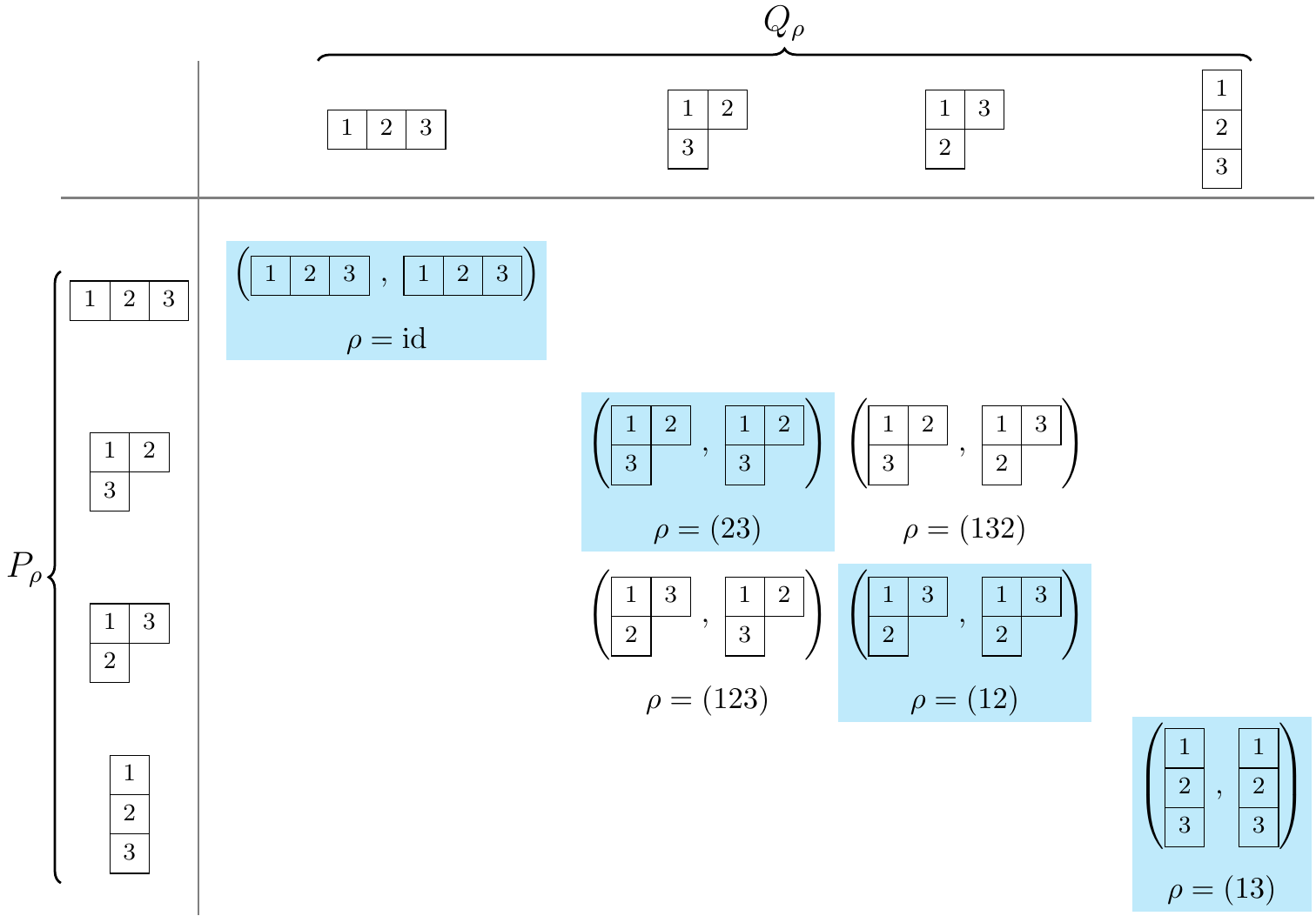}
  \caption{This table shows the elements in $\mathcal{Y}_3^{PQ}$, the
    set of all ordered pairs $(P_{\rho},Q_{\rho})$ of tableaux of the
    same shape, where $P_{\rho},Q_{\rho}\in\mathcal{Y}_3$. Each of the
    pairs $(P_{\rho},Q_{\rho})$ corresponds to a unique permutation
    $\rho\in S_3$ via the Robinson-Schensted correspondence; this
    permutation is captured together with the pair of tableaux
    here. The diagonal (shaded) elements are those that obey
    $P_{\rho}=Q_{\rho}$, which are clearly in $1$-to-$1$
    correspondence with the elements of $\mathcal{Y}_3$ (\emph{c.f.}
    eq.~\eqref{eq:Yk-isom-pairs-sam-tableaux}). Notice that all
    permutations corresponding to pairs $(P_{\rho},P_{\rho})$ on the
    diagonal are involutions, and none of the permutations on the
    off-diagonal are involutions, as is claimed in
    Theorem~\ref{thm:Number-Involutions-Irreps}. Clearly, the block
    structure one obtains from arranging the pairs
    $(P_{\rho},Q_{\rho})$ in a table is the same as that obtained from
    arranging Hermitian projection and unitary transition operators in
    a table
    (\emph{c.f.}~\cite{Alcock-Zeilinger:2016sxc,Alcock-Zeilinger:2016cva}).
    The proof of Theorem~\ref{thm:Number-Involutions-Irreps} shown
    below clarifies that the here exemplified correspondence between
    the pairs of tableaux in $\mathcal{Y}_k^{PQ}$ and the Hermitian
    projection and unitary transition operators in $\mathfrak{S}_k$ is
    true in general.}
  \label{fig:Tableaux-Pairs}
\end{figure}

The Robinson-Schensted (RS)
prescription~\cite{Schensted:1961,Robinson:1938} defines a bijection
between the symmetric group $S_k$ and the set $\mathcal{Y}_k^{PQ}$ of \emph{all} pairs
of Young tableaux \emph{of the same shape}
\begin{equation}
  S_k
  \stackrel{\text{RS}}{\cong}
  \mathcal{Y}_k^{PQ}
  :=
  \left\lbrace
    (P,Q) \in \mathcal{Y}_k \times \mathcal{Y}_k
    \middle\vert
    \text{$P$ and $Q$ have the same shape}
  \right\rbrace
\end{equation}
(for the definition of the Robinson-Schensted bijection refer to,
e.g.,~\cite{Schensted:1961,Robinson:1938,Sagan:2000}).
Ref.~\cite{Schensted:1961} refers to the unique pair of tableaux
corresponding to a particular permutation $\rho\in S_k$ as the
$P$-symbol and the $Q$-symbol of $\rho$; we will denote the unique
ordered pair of tableaux corresponding to $\rho$ by
$(P_{\rho},Q_{\rho})$ in order to make the permutation $\rho$
explicit. Notice that the subset of $\mathcal{Y}_k^{PQ}$ in which
$P_{\rho}=Q_{\rho}$ is isomorphic to $\mathcal{Y}_k$,
\begin{equation}
  \label{eq:Yk-isom-pairs-sam-tableaux}
  \mathcal{Y}_k
  \cong
  \left\lbrace
    (P,P)
    \middle\vert
    P \in \mathcal{Y}_k
  \right\rbrace
  \subset
  \mathcal{Y}_k^{PQ}
\end{equation}
(see Figure~\ref{fig:Tableaux-Pairs} for an example
illustrating the set $\mathcal{Y}_k^{PQ}$ and
eq.~\eqref{eq:Yk-isom-pairs-sam-tableaux}). It remains to show that
the permutations in $S_k$ that map to the set $\mathcal{Y}_k^{PP}$ are
exactly those that are their own inverse.

It is shown in~\cite{Schuetzenberger:1963} that, under the RS
prescription, the $P$-symbol of the inverse permutation $\rho^{-1}$ is
$Q_{\rho}$, and the $Q$-symbol of $\rho^{-1}$ is $P_{\rho}$, that is
\begin{equation}
  P_{\rho^{-1}} = Q_{\rho}
  \qquad \text{and} \qquad
  Q_{\rho^{-1}} = P_{\rho}
  \ ,
\end{equation}
such that
\begin{equation}
  \rho \xleftrightarrow{\text{RS}} (P_{\rho},Q_{\rho})
  \qquad \text{and} \qquad
  \rho^{-1} \xleftrightarrow{\text{RS}} (P_{\rho^{-1}},Q_{\rho^{-1}})
  =  (Q_{\rho},P_{\rho})
  \ .
\end{equation}
If $\rho$ is in \emph{involution}, that is $\rho=\rho^{-1}$, then it's
$P$-symbol and $Q$-symbol must be equal. Since the Robinson-Schensted prescription
mapping every permutation $\rho\in S_k$ to a unique pair of tableaux
$(P_{\rho},Q_{\rho})$ (where $P_{\rho}$ and $Q_{\rho}$ have the same
shape) is a bijective map to begin with, we have now obtained a bijection between the
involutions of $S_k$ and the Young tableaux in $\mathcal{Y}_k$. With
all these prerequisites fully established, one obtains the following Theorem:

\begin{theorem}[Involutions and irreducible representations of $\SUN$]
\label{thm:Number-Involutions-Irreps}
The number of irreducible representations of $\SUN$ on a product space
$\Pow{k}$, where $N\geq k$, is given by the number of involutions in the symmetric group
$S_k$.
\end{theorem}

We remark that, if one interprets the elements of $S_k$ as linear maps
on $\Pow{k}$ (\emph{c.f.} section~\ref{sec:Invariant-Theory}), one may
forego the requirement that $N\geq k$ in
Theorem~\ref{thm:Number-Involutions-Irreps} and simply state that the
number of irreducible representations of $\SUN$ on a product space
$\Pow{k}$ is given by the maximum number of \emph{linearly
  independent} involutions in the symmetric group $S_k$.

The full details of the argument underlying
Theorem~\ref{thm:Number-Involutions-Irreps},
including proofs, is available in Donald Knuth's book ``The Art of
Computer Programming, Volume 3''~\cite[section~5.1.4]{Knuth:1998}.
Despite its elegance, the proof is rather lengthy: it takes up
approximately $10$ pages. Knuth himself comments: \emph{``This
  connection between involutions and tableaux is by no means obvious,
  and there is probably no very simple way to prove it''}.  The
present paper demonstrates that Knuth might have been more
optimistic: We present a proof (leading to a more general statement,
Corollary~\ref{cor:Irreps-SUN-MixedSpace}) that is less than one page
long. Our argument exploits concepts of representation theory, in
particular the theory of invariants.

\section{Invariant theory of $\SUN$ and the projector
  basis}\label{sec:Invariant-Theory}

Consider a product representations of $\SUN$ constructed from its
fundamental representation on a given vector space $V$ with
$\text{dim}(V)=N$, whose action will simply be denoted by
$v\mapsto U v$ for all $U\in\SUN$ and $v\in V$ (note that $U$
denotes both an element of $\SUN$ as well as its representation on $V$). Choosing a basis
$\{e_{(i)} |i = 1, \ldots, \text{dim}(V)\}$ such that
$v = v^i e_{(i)}$, this becomes $v^i \mapsto \tensor{U}{^i_j} v^j$.
This immediately induces a product representation of $\SUN$ on
$\Pow{k}$ if one uses this basis of $V$ to induce a basis on $\Pow{k}$
so that a general element $\bm v\in\Pow{k}$ takes the form
$\bm v = v^{i_1\ldots i_k} e_{(i_1)}\otimes\cdots\otimes e_{(i_k)}$:
\begin{equation}
  \label{eq:SUN-Action}
  (U \bm{v})^{i_1 \ldots i_k}
  :=
  \tensor{U}{^{i_1}_{j_1}} 
  \cdots 
  \tensor{U}{^{i_k}_{j_k}} 
  v^{j_1\ldots j_k} e_{(i_1)}\otimes\cdots\otimes e_{(i_k)}
  \ ,
\end{equation} 
denoting the representation of $U\in\SUN$ on $\Pow{k}$ by $U$ as well.
The theory of invariants~\cite{Weyl:1946}\footnote{The theory of
  invariants is also described in many modern textbooks such
  as~\cite{Bourbaki7-9:2000,Fulton:2004,Cvitanovic:2008zz,Tung:1985na},
  just to name a few.} exploits a set of functions (known as
\emph{invariants}) that are invariant under the action of the
group. That is, $\sigma$ is said to be an invariant of $\SUN$ on $\Pow{k}$ if, for
every $U\in\SUN$,
\begin{equation}
  \label{eq:invariant-of-SUN}
  \sigma \circ U = U \circ \sigma
  \ ,
\end{equation}
where $\circ$ denotes the composition of linear maps, and $U$ in
eq.~\eqref{eq:invariant-of-SUN} is understood to denote  the representations of the
group element $U\in\SUN$ on $\Pow{k}$. 

Since all the factors in $\Pow{k}$ are identical, the
notion of permuting factors is a natural one that leads to a map
from the permutation group $S_k$ to the space of invertible linear maps on
$\Pow{k}$, $\mathsf{GL}\left(\Pow{k}\right)$. In particular, we define the
action of a permutation $\rho\in S_k$ on $\Pow{k}$ 
by\footnote{Permuting the basis vectors instead involves $\rho$:
  $v^{i_{\rho^{-1}(1)}\ldots i_{\rho^{-1}(k)}} e_{(i_1)} \otimes \cdots \otimes
  e_{(i_k)} = v^{i_1\ldots i_k} e_{(i_{\rho(1)})} \otimes \cdots
  \otimes e_{(i_{\rho(k)})}$.}
\begin{equation}
  \label{eq:perm-facs-def}
  (\rho \bm{v})^{i_1 \ldots i_k} 
  :=
  v^{i_{\rho^{-1}(1)}\ldots i_{\rho^{-1}(k)}} 
  e_{(i_1)}\otimes\cdots\otimes e_{(i_k)}
  \ .
\end{equation}
If $\text{dim}(V)=N\geq k$, the map
$S_k\to\mathsf{GL}\left(\Pow{k}\right)$ defined
through~\eqref{eq:perm-facs-def} is bijective and we will simply refer
to the image of $\rho\in S_k$ by $\rho$. In the following, we
therefore think of $S_k$ as a subgroup of
$\mathsf{GL}\left(\Pow{k}\right)$. Notice that the
requirement $N\geq k$ is necessary to ensure that the linear
maps in $S_k$ do not become linearly dependent.

From the definitions~\eqref{eq:SUN-Action}
and~\eqref{eq:perm-facs-def} one immediately infers that the product
representation commutes with all permutations on any
$\bm v\in\Pow{k}$, $U \circ\rho\circ\bm{v}=\rho\circ U\circ \bm{v}$,
implying that any such permutation $\rho\in S_k$ is an
\emph{invariant} of $\SUN$ in accordance with
eq.~\eqref{eq:invariant-of-SUN}. It can further be shown that these
permutations in fact span the space of all linear invariants of $\SUN$
on $\Pow{k}$ (\cite[Theorem~13.7]{Tung:1985na} proves this statement
for $\mathsf{GL}(N)$ on $\Pow{m}$, but the proof can easily be adapted
to $\SUN$). The permutations are thus
referred to as the \emph{primitive invariants} of $\SUN$ on
$\Pow{k}$. The full space of linear invariants (viewed  as
a subspace of the space of linear maps on
$\Pow{k}$, $\Lin{\Pow{k}}$) is then given by
\begin{equation}
  \label{eq:PI-def}
  \API{\SUN,\Pow{k}} 
  := 
  \Bigl\{ 
  \sum_{\rho\in S_k} \alpha_\rho \sigma \Big| 
  \alpha_\rho\in \mathbb{C}, \rho\in S_k \Bigr\}
  \subset
  \Lin{\Pow{k}}
  \ .
\end{equation}
In~\cite{Alcock-Zeilinger:2016cva}, we explained that there exists a
basis of the algebra of invariants $\API{\SUN,\Pow{k}}$ in terms of
Hermitian Young projection operators onto the irreducible
representations of $\SUN$~\cite{Alcock-Zeilinger:2016sxc} and the
unitary transition operators between equivalent irreducible
representations (transition operators are essentially intertwining
operators acting on the whole space~\cite{Alcock-Zeilinger:2016cva},
\emph{c.f.}
eqns.~\eqref{eq:Trans-Ops-Def}).\footnote{Refs. ~\cite{Alcock-Zeilinger:2016sxc}
  and~\cite{Alcock-Zeilinger:2016cva} give compact construction
  algorithms for Hermitian Young projection operators and unitary
  transition operators, respectively.} For the sake of completeness,
we provide a short argument here:

Let $N\geq k$. Consider a set of Hermitian projection operators onto
the irreducible representations of $\SUN$ on $\Pow{k}$, and denote
this set by $\mathfrak{P}_{k}$. Suppose that the operators in
$\mathfrak{P}_{k}$ are mutually \emph{transversal} (following the
nomenclature of~\cite{Keppeler:2012ih}), that is to say that the
images of any two Hermitian projectors $P_i,P_j\in\mathfrak{P}_{k}$
only intersect at $0$,
\begin{equation}
  \label{eq:Transversality-Def}
  P_i P_j = \delta_{ij}
  \ .
\end{equation}
Then, clearly, $P_i$ and $P_j$ are orthogonal under the scalar product
\begin{equation}
  \label{eq:Scalar-Product-Def}
  \left\langle A\middle\vert B\right\rangle:=\Tr{A^{\dagger}B}
  \ .
\end{equation}
\begin{sloppypar}
Using Schur's Lemma~\cite{Schur:1901}, we may complete the set
$\mathfrak{P}_{k}$ to a basis of $\API{\SUN,\Pow{k}}$ (we cite a
combination of the versions of Schur's Lemma given
in~\cite{Schur:1901,Sagan:2000,KosmannSchwarzbach:2000}):
\end{sloppypar}

\begin{lemma}[Schur's Lemma]\label{lemma:Schur}
  Let $\mathcal{V}_i$ and $\mathcal{V}_j$ be two irreducible
  $\mathsf{G}$-modules of a group $\mathsf{G}$. Let
  $T_{ij}:\mathcal{V}_j\to\mathcal{V}_i$ be a $\mathsf{G}$-homomorphism. Then 
  \begin{enumerate}
  \item\label{itm:Schur1} $T_{ij}$ is a $\mathsf{G}$-isomorphism if and only if
    $\mathcal{V}_i$ and $\mathcal{V}_j$ carry equivalent
    representations of $\mathsf{G}$, or
  \item\label{itm:Schur2} $T_{ij}$ is the zero map.
  \end{enumerate}
\end{lemma}

Since each projection operator $P_i:\mathcal{V}_i\to\mathcal{V}_i$
defines an irreducible $\SUN$-module $\mathcal{V}_i$ on $\Pow{k}$,
Lemma~\ref{lemma:Schur} ensures us that for each pair of Hermitian
projection operators $P_i:\mathcal{V}_i\to\mathcal{V}_i$ and
$P_j:\mathcal{V}_j\to\mathcal{V}_j$ corresponding to equivalent
representations there exists a pair of transition operators
$T_{ij}:\mathcal{V}_j\to\mathcal{V}_i$ and $T_{ji}$ satisfying
\begin{subequations}
  \label{eq:Trans-Ops-Def}
  \begin{align}
    T_{ij} & = T_{ji}^{\dagger} \\
    P_i T_{ij} & = T_{ij} = T_{ij} P_j \\
    T_{ij} T_{ji}^{\dagger} & = P_i
  \end{align}
\end{subequations}
(eqns.~\eqref{eq:Trans-Ops-Def} describe what it means for $T_{ij}$ to
be an $\SUN$-isomorphism).  We will denote the set of transition
operators by $\mathfrak{T}_k$.  Thus, the set of projection operators
$\mathfrak{P}_k$ can be augmented by $\mathfrak{T}_k$ to form a set of
linearly independent operators $\mathfrak{P}_k\cup\mathfrak{T}_k$,
which are mutually orthogonal under the scalar
product~\eqref{eq:Scalar-Product-Def}. From part~\ref{itm:Schur2} of
Schur's Lemma, we know that there exist no further operators that
are linearly independent from the operators in
$\mathfrak{P}_k\cup\mathfrak{T}_k$. Thus,
$\mathfrak{P}_k\cup\mathfrak{T}_k$ constitutes a basis of
$\API{\SUN,\Pow{k}}$; we refer to this basis as the \emph{projector
  basis},\footnote{In~\cite{Keppeler:2012ih} such a basis is also
  called a \emph{multiplet basis}.} and denote it by $\mathfrak{S}_k$,
\begin{equation}
  \label{eq:Projectorbasis-Proj-Trans}
  \mathfrak{S}_k 
  = 
  \mathfrak{P}_k
  \cup
  \mathfrak{T}_k
  \ .
\end{equation}
Note that the requirement $N\geq k$ ensures that none of the operators
in $\mathfrak{S}_k$ vanish.\footnote{For $N<k$, it is possible that
  elements of $\mathfrak{S}_k$ become zero. Such operators are
  referred to as \emph{dimensionally null} operators and are discussed
  in more detail
  in~\cite{Alcock-Zeilinger:2016bss,Alcock-Zeilinger:2016sxc,Alcock-Zeilinger:2016cva}.}

For $N\geq k$, the symmetric group $S_k$ and the set of all Hermitian
projection and unitary transition operators of $\SUN$ on $\Pow{k}$,
$\mathfrak{S}_k$, both constitute a basis of $\API{\SUN,\Pow{k}}$,
implying that these two sets must have the same size,
\begin{equation}
  \label{eq:API-bases-permutations-projectors}
  \left\vert S_k \right\vert
  =
  \left\vert \mathfrak{S}_k \right\vert
  \ .
\end{equation}

We are now in a position to present a compact proof of
Theorem~\ref{thm:Number-Involutions-Irreps}:

\emph{Proof of Theorem~\ref{thm:Number-Involutions-Irreps}: \quad}
Consider the elements of $S_k$ to be linear maps on $\Pow{k}$ as
given in~\eqref{eq:perm-facs-def}. These linear maps are unitary
with respect to the scalar product~\eqref{eq:Scalar-Product-Def}, that
is $\rho^{-1} = \rho^{\dagger}$ for all $\rho\in S_k$.

Let $n_T\left(S_k\right)$ denote the number of non-Hermitian elements
in $S_k$, and let $n_P\left(S_{k}\right)$ denote the number of
Hermitian elements in $S_k$. Since $S_{k}$ is a group and all its
elements are unitary, the elements of $S_{k}$ either are involutions
$\rho^{\dagger}=\rho^{-1}=\rho$, or occur in Hermitian conjugate pairs
$(\rho,\rho^{\dagger})$ , forcing $n_T\left(S_{k}\right)$ to be an
even number. For each pair $(\rho,\rho^{\dagger})$ (where $\rho$ is
not an involution) we can construct a Hermitian element $h_{\rho}$ and an
anti-Hermitian element $a_{\rho}$ as
  \begin{equation}
    \label{eq:HermitianConjPairs1}
    h_{\rho} := \rho + \rho^{\dagger} 
    \qquad \text{and} \qquad 
    a_{\rho} := \rho - \rho^{\dagger}
    \ .
  \end{equation}
  Clearly, the set of all Hermitian elements of $S_k$ together with
  the set of all $h_{\rho}$ and $a_{\rho}$ for all non-Hermitian
  $\rho\in S_k$ constitutes a basis of $\API{\SUN,\Pow{k}}$, which
  will be denoted by $\Tilde{S}_k$. In particular, $\Tilde{S}_{k}$ splits into two
  disjoint subsets:
\begin{subequations}
  \label{eq:Hmn-Amn-size}
  \begin{enumerate}
  \item the set $H_{k}\subset \Tilde{S}_{k}$ consisting of all Hermitian
    elements in $\Tilde{S}_{k}$ with size
    \begin{equation}
      \label{eq:Hmn-size}
      \vert H_{k} \vert
      =
      n_P\left(S_{k}\right)
      +
      \frac{n_T\left(S_{k}\right)}{2}
      \ ,
    \end{equation}
  \item the set $A_{k}\subset \Tilde{S}_{k}$ containing all anti-Hermitian
    elements in $\Tilde{S}_{k}$ with size
    \begin{equation}
      \label{eq:Amn-size}
      \vert A_{k} \vert
      =
      \frac{n_T\left(S_{k}\right)}{2}
      \ .
    \end{equation}
  \end{enumerate}
\end{subequations}
Similarly, the projector basis of $\API{\SUN,\Pow{k}}$,
$\mathfrak{S}_{k}$, consists of Hermitian elements (the projection
operators) and elements that occur in Hermitian conjugate pairs (the
transition operators). Thus, we may again form a set
$\Tilde{\mathfrak{S}}_{k}$ from $\mathfrak{S}_{k}$, such that
$\Tilde{\mathfrak{S}}_{k}$ is the union of two
disjoint sets $\mathfrak{H}_{k}$ and $\mathfrak{A}_{k}$ consisting of
Hermitian and anti-Hermitian elements, respectively, analogous to
eqns.~\eqref{eq:HermitianConjPairs1} and~\eqref{eq:Hmn-Amn-size}.
Since $H_{k}$ and $A_{k}$, respectively, $\mathfrak{H}_{k}$ and
$\mathfrak{A}_{k}$ are disjoint, it follows that
\begin{equation}
  \label{eq:HA-subsets-basis-invariants}
  \vert H_{k} \vert \overset{!}{=} \vert \mathfrak{H}_{k} \vert
  \qquad \text{and} \qquad
  \vert A_{k} \vert \overset{!}{=} \vert \mathfrak{A}_{k} \vert
\end{equation}
for both sets $\Tilde{S}_k=H_{k}\cup A_{k}$ and
$\Tilde{\mathfrak{S}}_{k}=\mathfrak{H}_{k}\cup\mathfrak{A}_{k}$ to constitute
a basis of the algebra of invariants $\API{\SUN,\Pow{k}}$.

Clearly, a linear combination of Hermitian objects will be Hermitian.
Thus, only the elements of $A_{k}$ and $\mathfrak{A}_{k}$ can be used
to construct the non-Hermitian elements of $\API{\SUN,\Pow{k}}$. From eq.~\eqref{eq:HA-subsets-basis-invariants} it
then follows that
\begin{equation}
\label{eq:HA-subsets-equal-size}
  n_T \left( \mathfrak{S}_{k} \right) =
  \vert \mathfrak{T}_{k} \vert =
  2 \vert \mathfrak{A}_{k} \vert =
  2 \vert A_{k} \vert =
  n_T \left( S_{k} \right)
   \ .
\end{equation}
Knowing that the sets $\mathfrak{S}_{k}$ and $S_{k}$ have the same
size (\emph{c.f.}  eq.~\eqref{eq:API-bases-permutations-projectors})
as do their subsets containing only the non-Hermitian elements
(eq.~\eqref{eq:HA-subsets-equal-size}), it follows that
  \begin{equation}
    \vert \mathfrak{P}_{k} \vert 
    =
    n_P\left(\mathfrak{S}_{k}\right)
    = 
    \vert \mathfrak{S}_{k} \vert - n_T\left(\mathfrak{S}_{k}\right)
    \xlongequal[\text{eq.~\eqref{eq:API-bases-permutations-projectors}}]{\text{eq.~\eqref{eq:HA-subsets-equal-size}}}
    \vert S_{k} \vert - n_T\left(S_{k}\right)
    = 
    n_P\left(S_{k}\right)
    \ ,
  \end{equation}
  thus concluding the proof of the theorem.\qed

\section{The irreducible representations
  of~\texorpdfstring{$\SUN$}{SU(N)} on mixed product spaces}\label{sec:Irreps-MixedSpace}

As explained in the previous section, the permutations in $S_k$ span
the algebra of invariants $\API{\SUN,\Pow{k}}$. These permutations may
be depicted graphically as
\emph{birdtracks}~\cite{Penrose:1971Com,Cvitanovic:2008zz}: Consider a
particular permutation $\rho\in S_k$. To obtain the birdtrack of
$\rho$, we write
two columns $(1,2,3,\ldots,k)^t$ next to each other, and then connect
the entry $i$ of the right column to the value of $\rho(i)$ in the
left column, marking each line with an arrow from right to left. We
then delete the numbers from the diagram, retaining only the
lines. For example,
\begin{equation}
  \rho = (134)(25)
  \hspace{2mm} \xlongrightarrow{\text{write columns}} \hspace{2mm}
  \begin{matrix}
    1 & 1 \\
    2 & 2 \\
    3 & 3 \\
    4 & 4 \\
    5 & 5
  \end{matrix}
  \hspace{2mm} 
  \xlongrightarrow{\text{draw lines}} 
  \hspace{2mm}
  \hspace{0.27mm}
  \raisebox{-8pt}{\FPic{5s134s25-MarkRho-Diagram}}
  \hspace{0.27mm}
  \hspace{2mm} 
  \xlongrightarrow[\text{only}]{\text{retain lines}} 
  \hspace{2mm}
  \FPic{5ArrLeft}
  \FPic{5s134s25SN}
  \FPic{5ArrRight}
  \ ;
\end{equation}
the last image is the birdtrack of $\rho$. Birdtracks ideally
lend themselves to be interpreted as linear maps on $\Pow{k}$
according to eq.~\eqref{eq:perm-facs-def}, for example the equation
\begin{equation}
  (123) \circ v_1\otimes v_2\otimes v_3 = v_3\otimes v_1\otimes v_2
\end{equation}
is written in the birdtrack formalism as
\begin{equation}
  \FPic{3ArrLeft}
  \FPic{3s123SN}
  \FPic{3ArrRight} \;
  \FPic{3v123Labels} 
  \; = \; 
  \FPic{3v312Labels} 
  \ ,
\end{equation}
where each term in the product $v_1\otimes v_2\otimes v_3$ (written as
a tower $\FPic{3v123Labels}$) can be thought of as being moved along
the lines of
$\FPic{3ArrLeft}\FPic{3s123SN}\FPic{3ArrRight}$~. Furthermore, the
birdtrack formalism allows for an efficient way to talk about
antifundamental representations of $\SUN$ on the mixed space
$\MixedPow{m}{n}$ and the associated algebra of invariants:

As one considers the fundamental representation of $\SUN$ on a vector
space $V$, one may also consider the \emph{anti-fundamental}
representation of $\SUN$ on the dual space $V^*$. Using the methods
outlined in section~\ref{sec:Invariant-Theory}, the irreducible
representations of $\SUN$ on a mixed product space $\MixedPow{m}{n}$
can again be determined using the invariants living in the
algebra~\cite{Cvitanovic:2008zz,Tung:1985na}
\begin{align}
  \label{eq:apimmp}
  \API{\SUN,\MixedPow{m}{n}} 
  & := \Bigl\{\sum_{\rho\in S_{m,n}} \alpha_\rho \rho \ \Big| 
  \ \alpha_\rho \in \mathbb{C} \Bigr\}
    \nonumber \\
  & \subset
  \Lin{\MixedPow{m}{n}}
  \ ,
\end{align}
where $S_{m,n}$ denotes the set of primitive invariants of $\SUN$ on
$\MixedPow{m}{n}$ as described below (\emph{c.f.}
eqns.~\eqref{eq:S3-to-S2plus1-map}). The elements of $S_{m,n}$ are in a
$1$-to-$1$ correspondence with the primitive invariants in $S_{m+n}$ as
we can construct them graphically from the elements of $S_{m+n}$ by
swapping the left and right endpoints on the specific $V$ in
$\Pow{(m+n)}$ to be converted into is dual vector space $V^*$. An
example will give clarity: The primitive invariants in $S_3$ map onto
those in $S_{2,1}$ as \allowdisplaybreaks[0]\begin{subequations}
\label{eq:S3-to-S2plus1-map}
\begin{IEEEeqnarray}{0lCCCCCCCCCCCr}
  S_3:\hspace{1cm}
  &
  \FPic{3ArrLeft}
  \FPic{3IdSN}
  \FPic{3ArrRight}
  &
  \ , \quad 
  &
  \FPic{3ArrLeft}
  \FPic{3s12SN}
  \FPic{3ArrRight}
  &
  \ , \quad 
  &
  \FPic{3ArrLeft}
  \FPic{3s23SN}
  \FPic{3ArrRight}
  &
  \ , \quad 
  &
  \FPic{3ArrLeft}
  \FPic{3s13SN}
  \FPic{3ArrRight}
  &
  \ , \quad
  &
  \FPic{3ArrLeft}
  \FPic{3s123SN}
  \FPic{3ArrRight}
  &
  \ , \quad
  &
  \FPic{3ArrLeft}
  \FPic{3s132SN}
  \FPic{3ArrRight}
  & 
  \\
  &
  \FPic{StraightDownArrow}
  & & 
  \FPic{StraightDownArrow}
  & & 
  \FPic{StraightDownArrow}
  & & 
  \FPic{StraightDownArrow}
  & & 
  \FPic{StraightDownArrow}
  & & 
  \FPic{StraightDownArrow}
  &
  \nonumber \\
  S_{2,1}:
  &
  \FPic{2q1qbArr}
  \FPic{3IdSN}
  \FPic{2q1qbArr}
  &
  \ , \quad 
  &
  \FPic{2q1qbArr}
  \FPic{3s12SN}
  \FPic{2q1qbArr}
  &
  \ , \quad 
  &
  \FPic{2q1qbArr}
  \FPic{2q1qbTl23r23N}
  \FPic{2q1qbArr}
  &
  \ , \quad 
  &
  \FPic{2q1qbArr}
  \FPic{2q1qbTl13r13N}
  \FPic{2q1qbArr}
  &
  \ , \quad
  &
  \FPic{2q1qbArr}
  \FPic{2q1qbTl13r23N}
  \FPic{2q1qbArr}
  &
  \ , \quad
  &
  \FPic{2q1qbArr}
  \FPic{2q1qbTl23r13N}
  \FPic{2q1qbArr}
  & 
  \label{eq:q2qb1}
  \ .
\end{IEEEeqnarray}
\end{subequations}%
\allowdisplaybreaks[1]%
Since the primitive invariants of $\SUN$ on $\MixedPow{m}{n}$ are
Hermitian with respect to the scalar
product~\eqref{eq:Scalar-Product-Def} if and only if their birdtrack
expressions are symmetric under a reflection about the vertical
axis\footnote{This statement is only true for the \emph{primitive}
  invariants: While any birdtrack operator that is symmetric under a
  flip about the vertical axis is Hermitian, the converse need not
  hold. In particular, the birdtrack formalism allows for a compact
  description of operators in the algebra of invariants, which, albeit
  being Hermitian, may not be explicitly symmetric when mirrowed about
  their vertical axis. However, once these operators are resolved into
  their constituent primitive invariants, their Hermiticity properties
  are once again uniquely determined through a reflection about the
  vertical
  axis~\cite{Alcock-Zeilinger:2016bss,Alcock-Zeilinger:2016sxc,Cvitanovic:2008zz}.}
(\emph{c.f.}~\cite{Cvitanovic:2008zz,Alcock-Zeilinger:2016sxc}), the
graphical procedure of transforming some of its fundamental legs into
antifundamental legs does not affect the Hermiticity of the
birdtrack. Thus, the subset of Hermitian elements in $S_{m,n}$ and
$S_{m+n}$ have the same size.

It should be noted that $S_{m,n}$, unlike $S_{m+n}$, is not a group;
for example, the last four elements of $S_{2,1}$ in~\eqref{eq:q2qb1} do
not have inverses.

\begin{sloppypar}
  The algebra of invariants $\API{\SUN,\MixedPow{m}{n}}$
  again allows for an orthogonal basis in terms of Hermitian
  projection and unitary transition operators, as the argument based
  on Schur's Lemma (Lemma~\ref{lemma:Schur}) given in
  section~\ref{sec:Invariant-Theory} can immediately be adapted to
  this situation: Lemma~\ref{lemma:Schur} tells us that the set of
  projection and transition operators constitute the largest possible
  set of linearly independent operators.  (Alternatively, an argument
  in terms of Clebsch-Gordan operators is given
  in~\cite{Alcock-Zeilinger:2017Singlets}). The projector basis of
  $\API{\SUN,\MixedPow{m}{n}}$ will be denoted by
  $\mathfrak{S}_{m,n}$.
\end{sloppypar}

Therefore, since $S_{m,n}$ and
$S_{m+n}$ contain the same number of Hermitian elements, and since
$\mathfrak{S}_{m,n}$ and $S_{m,n}$ each constitute a basis of the algebra of
invariants $\API{\SUN,\MixedPow{m}{n}}$, the proof of
Theorem~\ref{thm:Number-Involutions-Irreps} immediately translates to
the irreducible representations of $\SUN$ on the mixed space
$\MixedPow{m}{n}$. This gives rise to the following Corollary:

\begin{corollary}[Irreducible representations of $\SUN$ on
  $\MixedPow{m}{n}$]\label{cor:Irreps-SUN-MixedSpace}
  The number of irreducible representations of $\SUN$ over a product
  space $\MixedPow{m}{n}$ does not depend on $m$ and $n$ individually,
  but only on the sum $(m+n)$.
  In particular, the number of irreducible representations of $\SUN$ on any space
  $\MixedPow{m}{n}$ satisfying $m+n=k$ is the same as the number of
  irreducible representations of $\SUN$ on $\Pow{k}$.
\end{corollary}

\section{Conclusion}\label{sec:Conclusion}

In the present paper, we provided an alternative, simple proof of
Theorem~\ref{thm:Number-Involutions-Irreps}, which states that the
number of irreducible representations of $\SUN$ on $\Pow{k}$ is the
same as the number of involutions in the symmetric group $S_k$. The
proof given in section~\ref{sec:Invariant-Theory} heavily utilizes
the theory of invariants. In particular, we rely on the fact that the
set of Hermitian projection operators on the irreducible
representations of $\SUN$ on $\Pow{k}$, and the unitary transition
operators between equivalent representations, span the algebra of
invariants $\API{\SUN,\Pow{k}}$ (this last fact follows directly from
Schur's Lemma, but was als proven independently of Schur's Lemma 
in~\cite{Alcock-Zeilinger:2016cva}).

In section~\ref{sec:Irreps-MixedSpace}, we used the birdtrack
formalism to show that the number of Hermitian primitive invariants of
$\SUN$ on $\Pow{k}$ is the same as the number of Hermitian primitive
invariants on $\MixedPow{m}{n}$ for $m+n=k$ (we understand Hermiticity
with respect to the scalar
product~\eqref{eq:Scalar-Product-Def}). This, together with the fact
that the set of Hermitian projection and unitary transition operators
$\mathfrak{S}_{m,n}$ also constitutes a basis for
$\API{\SUN,\MixedPow{m}{n}}$, allowed us to formulate
Corollary~\ref{cor:Irreps-SUN-MixedSpace}, which states that the
number of irreducible representations of $\SUN$ on any product space
consisting of factors $V$ and $V^*$ only depends on the total number
of factors, but not on which of them are $V$ and which are its dual
$V^*$.

Corollary~\ref{cor:Irreps-SUN-MixedSpace} has interesting
consequences: As was explained in section~\ref{sec:Introduction}, the
irreducible representations of $\SUN$ on $\Pow{k}$ are classified by
Young tableaux. Similarly, one may construct tableaux corresponding to
the irreducible representations of $\SUN$ on $\MixedPow{m}{n}$ using
an algorithm that goes back to Littlewood and
Richardson~\cite{LittlewoodRichardson:1934}. Corollary~\ref{cor:Irreps-SUN-MixedSpace}
states that the Littlewood-Richardson construction must yield the same
number of tableaux on $\MixedPow{m}{n}$ and $\MixedPow{m'}{n'}$ if
$m+n=m'+n'$. In particular, the number of Littlewood-Richardson
tableaux of $\SUN$ on $\MixedPow{m}{n}$ must be the same number as
$\left\vert\mathcal{Y}_{m+n}\right\vert$, the number of Young tableaux
consisting of $m+n$ boxes.

\paragraph{Acknowledgements:} H.W. is supported by South Africa's
National Research Foundation under CPRR grant number 90509. J.A-Z. was
supported (in sequence) by the postgraduate funding office of the
University of Cape Town (2014), the National Research Foundation
(2015) and the Science Faculty PhD Fellowship of the University of
Cape Town (2016-2017).

 \bibliographystyle{utphys}
 \bibliography{GroupTheory,PaperLibrary,BookLibrary,MiscLibrary}

\end{document}